\documentclass[journal,twoside,web]{ieeecolor}

\usepackage{generic}
\usepackage{cite}
\usepackage{amsmath,amssymb,amsfonts}
\usepackage{algorithmic}
\usepackage{graphicx}
\usepackage{textcomp}
\usepackage{breqn}
\usepackage{graphicx}
\usepackage{subfigure}
\usepackage{epsfig} 
\usepackage{cite}
\usepackage{amsthm}  
\usepackage{lcsys}

\newtheorem{theorem}{Theorem}[section]
\newtheorem{lemma}[theorem]{Lemma}

\newtheorem{definition}[theorem]{Definition}
\newtheorem{remark}[theorem]{Remark}
\newtheorem{assumption}[theorem]{Assumption}
{\theoremstyle{plain}
\newtheorem*{dem}{Proof}}
\newcommand{\RNum}[1]{\uppercase\expandafter{\romannumeral #1\relax}}
\newcommand{\normm}[1]{{\left\vert\kern-0.25ex\left\vert\kern-0.25ex\left\vert #1
    \right\vert\kern-0.25ex\right\vert\kern-0.25ex\right\vert}}

\begin{document}

\def\BibTeX{{\rm B\kern-.05em{\sc i\kern-.025em b}\kern-.08em
    T\kern-.1667em\lower.7ex\hbox{E}\kern-.125emX}}
\markboth{\journalname, VOL. XX, NO. XX, XXXX 2017}
{Author \MakeLowercase{\textit{et al.}}: Preparation of Papers for IEEE Control Systems Letters (August 2022)}

\title{Method of Successive Approximations for Stochastic Optimal Control: Contractivity and Convergence}

\author{
Safouane TAOUFIK$^{a}$,
Badr MISSAOUI$^{b}$
\thanks{$^{a}$ College of Computing, UM6P; safouane.taoufik@um6p.ma}
\thanks{$^{b}$ Moroccan Center of Game Theory, UM6P;badr.missaoui@um6p.ma}
}

\maketitle
\thispagestyle{empty}
\begin{abstract}
The Method of Successive Approximations (MSA) is a fixed-point iterative method used to solve stochastic optimal control problems. It is an indirect method based on the conditions derived from the Stochastic Maximum Principle (SMP), an extension of the Pontryagin Maximum Principle (PMP) to stochastic control problems. In this study, we investigate the contractivity and the convergence of MSA for a specific and interesting class of stochastic dynamical systems (when the drift coefficient is one-sided-Lipschitz with a negative constant and the diffusion coefficient is Lipschitz continuous). Our analysis unfolds in three key steps: firstly, we prove the stability of the state process with respect to the control process. Secondly, we establish the stability of the adjoint process. Finally, we present rigorous evidence to prove the contractivity and then the convergence of MSA. This study contributes to enhancing the understanding of MSA's applicability and effectiveness in addressing stochastic optimal control problems. 
\end{abstract}


\section{Introduction}
\label{sec:introduction}
\IEEEPARstart{S}{tochastic} optimal control problems find a wide range of applications, including finance, resource management, robotics, and autonomous systems. However, their analytical solutions are frequently extremely challenging and, in many cases, nearly impossible to obtain, except for very specific cases. Therefore, the need for numerical methods becomes apparent. Among these numerical methods, the Method of Successive Approximation (MSA) emerges as an effective method for tackling optimal control challenges. MSA is an iterative approach for solving stochastic control problems and is rooted in Pontryagin's optimality principle.

Recent works have expanded the usage of the MSA to address a deep learning problem. This presents a novel approach for training deep neural networks by applying the principles of optimal control \cite{weinanE}.

The paper seeks to prove the convergence of MSA towards the solution of the stochastic optimal control problem. To do so, we build upon the work of Smith and Bullo \cite{smith}, which focused on studying the convergence of the MSA's in the deterministic optimal control. Our contribution extends this work to address the case of the stochastic control problem. 

Our contributions are as follows: we demonstrate, under some assumptions (to be clearly specified later), that the function mapping a given control input to the system's response (the associated state process) is Lipschitz continuous with respect to a norm that will be explicitly defined. Additionally, we establish the same property for adjoint process, which appears in the Stochastic Maximum Principle, by using its boundedness under specific assumptions that will also be clarified. Finally, based on the results mentioned earlier, we establish that the MSA function is a contraction map in a closed set of a Banach space.
\section{Preliminaries}
\label{sec:preliminaries}

\subsection{Notations:}
\noindent \textbf{For vectors:} We use canonical scalar product and norm.
\\For $x,y\in \mathbb{R}^n$: $\langle x,y \rangle:=x^{T}y$, $\|x\|:=\sqrt{ \langle x,x \rangle}$.
\\\textbf{For matrices:} 
For $A\in \mathcal{M}_{n,m}(\mathbb{R})$, we denote the induced 2-norm by $\normm{A}$ and the Frobenius norm by $\|A\|:=\sqrt{\operatorname{Tr}(A^TA)}$.\\
For $A\in\mathcal{M}_{n} (\mathbb{R})$, we denote by $\mu(A)$ the induced logarithmic norm, defined as: $\mu(A):=\displaystyle\lim_{\alpha \to 0^{+}} \frac{\normm{I_n+\alpha A}-1}{\alpha}$.
\\\textbf{For functions:} For $f:\mathbb{R} \to \mathbb{R}$, we use the notations $D^{+}f(t)$, $D_{+}f(t)$, $D^{-}f(t)$ and $D_{-}f(t)$ to represent the upper right, lower right, upper left, and lower left Dini derivatives, respectively.
\subsection{Some tools:}
\noindent Here are some Known results about Dini derivatives, the induced logarithmic norm and the one-sided-Lipschitz functions that we will use. (refer to \cite{MVT} and \cite{bullobook}) 
\noindent
 \begin{theorem} \textbf{(Dini derivatives and local strict monotonicity)}

\begin{enumerate}
    \item if $D_-f(t)>0$ then $f$ is locally strictly increasing to the left of $t$, i.e there exists $\eta_1>0$ such that:
    $$t-\eta_1<x<t \implies f(x)<f(t). $$
    \item if $D_+f(t)>0$ then $f$ is locally strictly increasing to the right of $t$, i.e there exists $\eta_2>0$ such that: $$t<x<t+\eta_2 \implies f(t)<f(x). $$
\end{enumerate}

\end{theorem}

 \begin{theorem} \textbf{(Mean Value Theorem for Dini Derivatives)}
Let $f:\mathbb{R}\rightarrow \mathbb{R}$ be a continuous function and $a<b$, then 
every Dini derivative upper bounds the quantity $$\gamma:=\displaystyle\frac{f(b)-f(a)}{b-a},$$ at one point and lower bounds it at another point,
which is equivalent to the existence of $ c_1,c_2,c_3,c_4\in ]a,b[,$ such that:
$$D^{+}f(c_1)\leq \gamma,D^{-}f(c_2)\leq \gamma,\gamma \leq D_{+}f(c_3), \gamma\leq D_{-}f(c_4).$$ 
\end{theorem}

\begin{theorem} \textbf{(Dini derivatives and monotonicity)}
\\Let $f$ be a continuous function.
\begin{enumerate}
    \item If one Dini derivatives is always non-negative, then $f$ is non-decreasing.
    \item If one Dini derivatives is always non-positive, then $f$ is non-increasing.
\end{enumerate}
\end{theorem}
 \begin{theorem} \textbf{(Gronwall inequality and Dini derivatives)}
 \\If $f$,  $a$ and $g$ are continuous functions such that:\\$\forall t, D^{+} f(t)\leq a(t)f(t)+g(t)$ then $$\forall t\geq t_0, f(t)\leq \left[ f(t_0)+\int_{t_0}^{t} g(s)e^{-A(s)}ds\right]e^{A(t)},$$ where $A$ is the antiderivative of $a$ that vanishes at $t_0$.
\end{theorem}
\noindent
\begin{theorem}
\textbf{(The induced logarithmic norm and inner product)}
\newline
Let $A\in \mathcal{M}_{n}(\mathbb{R})$. We have:

$$\mu(A)=\displaystyle\sup_{x\neq 0}\frac{\langle Ax,x\rangle }{\|x\|^2 }.$$ 
\end{theorem}

\noindent
\begin{definition} \textbf{(One sided Lipschitz)}
\\Let $b:\mathbb{R}^{n}\rightarrow \mathbb{R}^{n}$. We say that $b$ is one-sided-Lipschitz if there exists $c \in \mathbb{R}$ such that for all $x$, we have: 
$$  \langle x-y , b(x)-b(y)  \rangle  \leq c \|x-y\|^{2}.$$

\end{definition} 

\noindent

 \begin{theorem}
\textbf{(one-sided-Lipschitz and the induced logarithmic norm.)}
\\Let $U\subset\mathbb{R}^{n}$ be an open convex set, and $b:U\rightarrow \mathbb{R}^{n}$ a function of class $\mathcal{C}^{1}$. 
Then : $$\sup_{x \neq y} \frac{ \langle x-y , b(x)-b(y)  \rangle} {\|x-y\|^{2}}=\sup_{x\in \mathbb{R}^{n}} \mu(D_xb),$$
where $D_xb$ is the Jacobian matrix of $b$ at $x$.
\end{theorem}

\subsection{Stochastic Optimal control problem:}
\noindent Let $(\Omega,\mathcal{P},\mathcal{F},\mathbb{F})$ be a filtered probability space.
We assume that the filtration $\mathbb{F}=(\mathcal{F}_t)_{t\in[0,T]}$ satisfies the usual conditions. Let $(W_t)_t$ be a $d-$dimensional Brownian motion defined on $(\Omega,\mathcal{P},\mathcal{F},\mathbb{F})$.
We consider the following state equation:
\begin{equation}dX_t=b(t,X_t,\alpha_t) dt + \sigma(t,X_t) dW_t, X_0=x_0\in \mathbb{R}^{n} ,\label{state}\end{equation}
and the following cost function 
\begin{equation} J(\alpha)=\mathbb{E}\left[\displaystyle\int_{0}^{T}\varphi(s,X_s,\alpha_s)\mathrm{d}s +\psi(X_T) \right], \label{cost}\end{equation}
where $X$ is the solution of  \eqref{state}, $\psi: \mathbb{R}^n\rightarrow \mathbb{R}$,
$\varphi: [0,T]\times\mathbb{R}^n\times A\rightarrow \mathbb{R}$,
$b: [0,T]\times\mathbb{R}^n\times A\rightarrow \mathbb{R}^n$, and 
$\sigma: [0,T]\times\mathbb{R}^n\rightarrow \mathcal{M}_{n,d}(\mathbb{R})$ are continuous functions. The action space $A$ is a compact set of $\mathbb{R}^{m}$. The space of admissible controls $\mathcal{A}$ that we consider is the set of  $A$-valued adapted processes. We adopt the notations $$\|\alpha_t\|_{L^2}:=\sqrt{\mathbb{E}[\|\alpha_t\|^2]}\text{ and }\|\alpha\|_{\mathcal{A}}:=\sup_{t\in [0,T]}\sqrt{\mathbb{E}[\|\alpha_t\|^2]}.$$
We suppose that for every admissible control $\alpha$ and initial state $x_0$ the equation \eqref{state} has a unique solution. For some examples of sufficient conditions on existence and uniqueness, we refer the reader to \cite{Yong, Nisio, carmona}. Note that besides existence and uniqueness of the solution, those conditions imply also that 
$\displaystyle \mathbb{E}[ \sup_{0 \leq t\leq T}\|X_t\|^2]<\infty$, which implies by dominated convergence theorem that the function $t\mapsto \mathbb{E}[\|X_t\|^2]$ is continuous.
\\Let's consider the Hamiltonian: \\$H:[0,T]\times\mathbb{R}^n\times\mathbb{R}^n\times\mathcal{M}_{n,d}(\mathbb{R})\times A\rightarrow \mathbb{R}^{n}$ defined by:
$$H(t,x,y,z,a)=y^Tb(t,x,a)+\operatorname{Tr}(\sigma(t,x)^Tz)+\varphi(t,x,a).$$
\noindent
\textbf{Stochastic maximum principle (SMP)} 

It is known (\cite{stochastic}, \cite{SMP}) that if $\alpha$ is an optimal control and $X$ the associated state process, solution of (1), then $\alpha$ satisfies the SMP, that is:
$$\forall t, \alpha_t\in \arg \min_{a\in A} H(t,X_t,Y_t,Z_t,a),$$
where $(Y,Z)$ is the solution to the associated adjoint equation \begin{equation} \label{adjoint} \begin{split}dY_t&=-D_xH(t,X_t,Y_t,Z_t,\alpha_t)dt + Z_tdW_t \\Y_T&=D_x\psi(X_T),\end{split}\end{equation}
which is equivalent to:
\begin{equation}
dY_t=[A_t+B_tY_t+C_t(Z_t)]dt + Z_tdW_t, Y_T=\xi, \label{adjoint with new notation}
\end{equation}
where
\begin{equation*}
\begin{split}
A_t:=-D_x\varphi(t,X_t,\alpha_t), B_t:=-D_xb(t,X_t,\alpha_t)^{T},\\C_t(Z_t):=-D_x(\mathrm{Tr}(\sigma(t,X_t)^{T}Z_t)), \xi:=D_x\psi(X_T).         
\end{split}
\end{equation*}
\subsection{The method of successive approximations (MSA):}
MSA is a fixed-point method used to find a control that satisfies the SMP, providing a good candidate for the solution of the stochastic optimal problem.\\\\
\fbox{%
\begin{minipage} {0.97\linewidth} \textbf{MSA}
\noindent \\Make an initial guess $\alpha^0$.
\\for $i=1,2,...,N$ 
\begin{enumerate}
    \item Find $X^i$ the state process of the control $\alpha^{i-1}$
    \item Find $(Y^i,Z^i)$ the solution to the adjoint equation associated to $X^{i}$ and $\alpha^{i-1}$ .
    \item Find $\alpha^i$ that satisfies\\ $\forall t, \alpha_t^i\in \arg \displaystyle\min_{a\in A} H(t,X^i_t,Y^i_t,Z^i_t,a).$
\end{enumerate}
return $\alpha^N$
\end{minipage}}\\\\
The goal of this work is to prove that this algorithm converges to the optimal control under some assumptions. To be more precise, we will prove that the function $\operatorname{MSA}:\mathcal{A}\rightarrow\mathcal{A}$ is a contraction,
where
$\operatorname{MSA}(\alpha)=\beta$ where $\beta$ is defined by: $\beta_t\in \arg \displaystyle\min_{a\in A} H(t,X_t,Y_t,Z_t,a)$ where $X$ is the state process associated to $\alpha$ and $(Y,Z)$ the solution of the associated adjoint equation.
\section{Stability of the state process:}
\noindent
\begin{assumption}
\noindent
\begin{enumerate}
    \item There exists $c\in\mathbb{R}$, 
  for all $t,x,\overline{x},\alpha$  $$\left\langle x-\overline{x},b(t,x,\alpha)-b(t,\overline{x},\alpha)\right \rangle \leq c \|x-\overline{x} \|^2. $$
    \item There exists $L^{\alpha}_{b} \geq 0$, for all $ t,x,\alpha,\overline{\alpha}$  $$\|b(t,x,\alpha)-b(t,x,\overline{\alpha})\|\leq L^{\alpha}_{b} \|\alpha-\overline{\alpha} \|.$$ 
    \item There exists $L^{x}_{\sigma} \geq 0$, for all $t,x,\overline{x}$ $$ \|\sigma(t,x)-\sigma(t,\overline{x})\|\leq L^{x}_{\sigma} \|x-\overline{x}\|.$$ 
    \item $\mu:=-(c+\displaystyle\frac{(L^{x}_{\sigma})^2}{2}) > 0.$
\end{enumerate}
\end{assumption}
\begin{remark}
\noindent\\
Despite the fact that the last assumption ($\mu>0$) may seem unusual, it has already been used in the study of stability of numerical methods for systems of stochastic differential equations (refer to \cite{schurz} and \cite{meansquare}).
\end{remark}
\noindent
\begin{theorem}
Let $X$ and $\overline{X}$ be two solutions of \eqref{state} corresponding to controls $\alpha$ and $\overline{\alpha}$, and to the initial conditions $X_0$ and $\overline{X}_0$ respectively. Under the previous assumptions, we have:
\begin{equation} \label{state stability}
\begin{split}\|X_t-\overline{X}_t\|_{L^2}\leq & \|X_0-\overline{X}_0\|_{L^2}e^{-\mu t}\\&+L^{\alpha}_{b}\frac{1-e^{-\mu t}}{\mu} \|\alpha-\overline{\alpha}\|_{\mathcal{A}}.
\end{split}
\end{equation}
\end{theorem}
\noindent Before we proceed with the proof of Theorem \RNum{3}.3, we will establish the following lemma, that will not only be used in this proof but will also be instrumental in proving another theorem later on.
\begin{lemma} 
If $f$ and $g$ are continuous and non-negative functions such that for every $t$ that is not a zero of $f$, we have $D^+f(t)\leq g(t)$, then  this inequality holds for every $t$.

\end{lemma}

\begin{dem}  \textbf{(Lemma III.4)}
\\Let's first prove that $\forall t, D_+f(t)\leq g(t)$.
\\If $f(t)\neq 0$, there is nothing to prove.
\\If $f(t)=0$, we have $0\leq \displaystyle\frac{f(t+h)}{h}$, which implies \\$0\leq D_+f(t)$.\\
\noindent So, we have two possible cases $D_+f(t)=0$ or $D_+f(t)>0$. \\
For the first case, since $g(t)\geq 0$, then $D_+f(t)\leq g(t)$. \\
For the second case, $f$ is locally strictly increasing to the right of $t$ (i.e for every $h$ sufficiently small, $0<f(t+h)$ ), which implies that there is no zero of $f$ near $t$ at the right. By the Mean Value Theorem for Dini Derivatives, for any $h$, there exists $c_h\in ]t,t+h[$ such that $$\frac{f(t+h)}{h}\leq D_+f(c_h).$$ Since $c_h$ is not a zero of $f$ (for $h$ sufficiently small), we have $$\frac{f(t+h)}{h}\leq D_+f(c_h)\leq g(c_h),$$ and by taking the limit inferior as $h$ approaches 0 (due to the continuity of $g$), we obtain $D_+f(t)\leq g(t)$.\\
In both cases, we have $D_+f(t)\leq g(t)$, so we have for all $t$, $D_+f(t)\leq g(t)$.
\\\noindent Now, let's prove the result. For $t$  such that $f(t)=0$, and $h$ sufficiently small, by Mean Value Theorem, there exists $c_h\in ]t,t+h[$ such that $$\frac{f(t+h)}{h}\leq D_+f(c_h).$$ Based on the previous result, regardless of whether $c_h$ is a zero of $f$ or not, we always have: 
$$\frac{f(t+h)}{h}\leq D_+f(c_h)\leq g(c_h).$$
As we take the limit superior as $h$ approaches 0 (due to the continuity of $g$), we obtain the result.

\end{dem}
\begin{dem}\textbf{(Theorem III.3)} 
Let $X$ and $\overline{X}$ be two solutions of \eqref{state} corresponding to controls $\alpha$ and $\overline{\alpha}$, respectively. 
\\Applying Ito's formula to $\|X_t-\overline{X}_t\|^2$, and taking the expectation, we find for $s<t$ 
\begin{equation*}
\begin{split}
\mathbb{E}[\|X_t&-\overline{X}_t\|^2]-\mathbb{E}[\|X_s-\overline{X}_s\|^2]\\=&\mathbb{E}[ \int_s^t (2 \langle X_r-\overline{X}_r,b(r,X_r,\alpha_r)-b(r,\overline{X}_r,\overline{\alpha}_r) \rangle \\&+\|\sigma(r,X_r)-\sigma(r,\overline{X}_r)\|^2) dr],\end{split} \end{equation*}
According to the assumptions, we have \begin{equation*}
\|\sigma(r,X_r)-\sigma(r,\overline{X}_r)\|^2\leq (L^{x}_{\sigma})^2 \|X_r-\overline{X}_r\|^2.\end{equation*}
We also have
\begin{equation*}
\begin{split}
\langle X_r-&\overline{X}_r,b(r,X_r,\alpha_r)-b(r,\overline{X}_r,\overline{\alpha}_r)\rangle\\&=\langle X_r-\overline{X}_r,b(r,X_r,\alpha_r)-b(r,\overline{X}_r,\alpha_r) \rangle\\& + \langle X_r-\overline{X}_r,b(r,\overline{X}_r,\alpha_r)-b(r,\overline{X}_r,\overline{\alpha}_r)\rangle,  
\end{split}
\end{equation*} 
Since $b$ is one-sided-Lipschitz, therefore $$ \langle X_r-\overline{X}_r,b(r,X_r,\alpha_r)-b(r,\overline{X}_r,\alpha_r) \rangle \leq c \|X_r-\overline{X}_r \|^2.$$
By Cauchy-Schwarz inequality, and by the fact that $b$ is Lipschitz with respect to $\alpha$:
\begin{equation*}
\begin{split}
\langle X_r-&\overline{X}_r,b(r,\overline{X}_r,\alpha_r)-b(r,\overline{X}_r,\overline{\alpha}_r)\rangle \\& \leq\|X_r-\overline{X}_r\| \|b(r,\overline{X}_r,\alpha_r)-b(r,\overline{X}_r,\overline{\alpha}_r)\|\\ &\leq L^{\alpha}_{b} \|X_r-\overline{X}_r\| \|\alpha_r-\overline{\alpha}_r\|.\end{split}
\end{equation*} 
\noindent Combining all of the above, we obtain
\begin{equation*}
\begin{split}   
&\mathbb{E}[\|X_t-\overline{X}_t\|^2]-\mathbb{E}[\|X_s-\overline{X}_s\|^2]\leq\\&\mathbb{E}[\displaystyle\int_s^t (-2\mu \|X_r-\overline{X}_r\|^2+2 L^{\alpha}_{b} \|X_r-\overline{X}_r\| \|\alpha_r-\overline{\alpha}_r\|)dr].
\end{split} 
\end{equation*}
We use \begin{equation*}L_b^{\alpha}\mathbb{E}[\|X_r-\overline{X}_r\| \|\alpha_r-\overline{\alpha}_r\|]\leq f(r)g,
\end{equation*}
where
\begin{equation*}
f(t):=\|X_t-\overline{X}_t\|_{L^2} \text{, }g:=L^{\alpha}_{b}\|\alpha-\overline{\alpha}\|_{\mathcal{A}}.
\end{equation*}
Dividing by $t-s$ and taking the limit superior as $t$ approaches $s^+$, we obtain $$D^+ f(t)^2\leq -2\mu f(t)^2+2f(t)g.$$
Let's prove that for every $t$, \begin{equation}\label{ineq Dini}
    D^+ f(t)\leq -\mu f(t)+g.
\end{equation} \\If $f(t)\neq 0$ (so $f(t)>0$), then
$$\frac{f(t+h)^2-f(t)^2}{h}=(f(t+h)+f(t))\frac{f(t+h)-f(t)}{h}.$$
Since $f$ is continuous at $t$, we have  $$\displaystyle\lim_{h\to 0^+} f(t+h)+f(t)=2f(t)>0,$$ and then $$D^{+} f(t)^2= 2 f(t)D^{+}f(t).$$ hence the result.
\\Now, if $f(t)=0$, then \eqref{ineq Dini} is equivalent to $D^+ f(t)\leq g,$ which holds also when $f(t)\neq 0$ (by the previous case and by the fact that $-\mu<0$), so we conclude by the lemma.\\
\noindent By Gronwall's lemma (Theorem \RNum{2}.4), we have that $$f(t)\leq f(0) e^{-\mu t}+\int_0^t g e^{-\mu(t-s)} ds.$$
Substituting the expressions of $f$ and $g$, we obtain the result.
\end{dem}

\section{Boundedness of the adjoint process}
\noindent Alongside the prior assumptions, we introduce the following additional assumptions for this section:
\begin{assumption}
\noindent
\begin{enumerate}
\item There exists $M \geq 0$ such that $\|D_x\psi\| \leq M$.
\item There exists $a \geq 0$ such that $\|D_x\varphi\| \leq a$. 
\end{enumerate}
\end{assumption}

\noindent We start by proving the following lemma:
\noindent
\begin{lemma}

\noindent For all $t,X,Z$, we have
$$\|D_x \operatorname{Tr}(\sigma(t,X)^T Z)\|\leq L^{x}_{\sigma} \|Z\|.$$
\end{lemma}
\begin{dem}
On one hand, if we expand the inequality, we get $$\sqrt{\sum_{k=1}^{n}(\sum_{l\in I}(\partial_k\sigma_{l}) Z_{l})^2}\leq L_{\sigma}^x\sqrt{\sum_{l\in I}Z_{l}^2},$$
where $I:=\{(i,j)\in[1..n]\times[1..d]\}$  
\\The inequality is equivalent to the 2-induced norm of the matrix $A=(\partial_k\sigma_{l})_{k,l}$  being less than $L^{x}_{\sigma}$.
\\On the other hand, if we expand the inequality resulting from the $L^{x}_{\sigma}$ Lipschitz continuity of $\sigma$ (for more precision, we take $y=x+\epsilon h$, divide the inequality by $\epsilon$, and then, as we let $\epsilon$ approach $0^+$), we obtain
$$\sqrt{\sum_{l\in I}(\sum_{k=1}^{n}(\partial_k\sigma_{l})h_{k})^2}\leq L_{\sigma}^x\sqrt{\sum_{k=1}^{n}h_{k}^2}.$$
This inequality is equivalent to the induced 2-norm of $A^T$ being less than $L^{x}_{\sigma}$.
We conclude by the fact that the induced 2-norm of a matrix is equal to the induced 2-norm of its transpose. 
\end{dem}

\begin{theorem}

Let $(Y,Z)$ be the solution of \eqref{adjoint with new notation}. Under the assumptions of this section, we have:

$$\|Y_t\|\leq \sqrt{e^{-\mu(T-t)}(M^2-(\frac{a}{\mu})^2)+(\frac{a}{\mu})^2}.
$$

\end{theorem}
\begin{dem}

\noindent 
\\By Cauchy-Schwarz inequality and by boundedness of $D_x\varphi$, we have:  \begin{equation} \label{A}
2\langle A_t,Y_t \rangle \geq -2a\|Y_t\|\geq -\frac{a^2}{\mu}-\mu\|Y_t\|^2.
\end{equation}
\\
By one-sided-Lipschitz condition of $b$ and using theorems \RNum{2}.5 and \RNum{2}.7, we have:
\begin{equation}
\langle B_tY_t,Y_t \rangle \geq -c \|Y_t\|^2. \label{Jacobian one sided}
\end{equation} 
\\By Cauchy-Schwarz inequality and by lemma \RNum{4}.2, we have:  \begin{equation} \label{C}
\langle C_t(Z_t),Y_t \rangle \geq -\frac{1}{2}\left((L^{x}_{\sigma})^2\|Y_t\|^2+\|Z_t\|^2\right).
\end{equation}
\\We apply Ito's formula to $e^{-\mu t}(\|Y_t\|^2-(\frac{a}{\mu})^2)$, and using  \eqref{A},\eqref{Jacobian one sided} and \eqref{C}, we obtain
\begin{equation} \label{inequality}
e^{-\mu T}(\|Y_T\|^2-(\frac{a}{\mu})^2)-e^{-\mu t}(\|Y_t\|^2-(\frac{a}{\mu})^2)\geq \int_t^T dM_s,\end{equation} where $$M_t:=\int_0^t 2e^{-\mu s} \langle Y_s,Z_sdW_s \rangle.$$
We have \begin{equation*}
    \|Z_s^TY_s\|\leq \normm{Z_s^T}\|Y_s\|\leq \|Z_s\|\|Y_s\|.
\end{equation*}
So
\begin{equation*} \begin{split}
\mathbb{E}[\sqrt{\int_0^T \|Z_s^TY_s\|^2 ds}] &\leq \sqrt{\mathbb{E}[\sup_t \|Y_t\|^2]} \sqrt{\mathbb{E}[\int_0^T \|Z_s\|^2 ds]}\\&<\infty.
\end{split}
\end{equation*}
And by Burkholder–Davis–Gundy inequality
\begin{equation*} 
\mathbb{E}[\sup_t |M_t| ]\leq C \mathbb{E}[\sqrt{\int_0^T \|Z_s^TY_s\|^2 ds}]<\infty.
\end{equation*}
So by dominated convergence theorem $(M_t)_t$ is a martingale.
\\Taking the conditional expectation with respect to $\mathcal{F}_t$ in \eqref{inequality}, we obtain: $$e^{-\mu t}(\|Y_t\|^2-(\frac{a}{\mu})^2)\leq \mathbb{E}[e^{-\mu T} (\|Y_T\|^2-(\frac{a}{\mu})^2))|\mathcal{F}_t].$$ 
By boundedness of $D_x\psi$, we obtain:
$$\|Y_t\|\leq \sqrt{e^{-\mu(T-t)}(M^2-(\frac{a}{\mu})^2)+(\frac{a}{\mu})^2}.
$$   
\end{dem}

\begin{remark}
This last quantity is clearly bounded \\(for example by:$\displaystyle \sqrt{M^2+(\frac{a}{\mu})^2}$ which does not depend on $T$).
\end{remark}
\section{Stability of the adjoint process:}
In addition to the previous assumptions, here we introduce the following additional assumptions for this section:

\begin{assumption} \noindent
\begin{enumerate}
    \item $D_xb$ and $D_x\varphi$ are Lipschitz with respect to $x$ and $\alpha$.
    \item $D_x\psi$ is Lipschitz with respect to $x$.
    \item $D_x\sigma$ does not depend on $x$.
\end{enumerate}
\end{assumption}
Let $X,Y,Z$ and $\overline{X},\overline{Y},\overline{Z}$ be two solutions of \eqref{state} and \eqref{adjoint} corresponding to the controls $\alpha$ and $\overline{\alpha}$ respectively.

\noindent
\begin{theorem}
We have :
\begin{equation}\label{adjoint stability}
\begin{split}    
\|&Y_t-\overline{Y}_t\|_{L^{2}}\leq L^{x}_{D_x\psi} \|X_T-\overline{X}_T\|_{L^{2}} e^{-\mu (T-t)}+\\&L_Y \displaystyle\int_t^T (\|X_{T+t-s}-\overline{X}_{T+t-s}\|_{L^{2}}+\|\alpha-\overline{\alpha}\|_{\mathcal{A}}) e^{-\mu (T-s)} ds.
\end{split}
\end{equation}
\end{theorem}
\begin{dem}
We apply Ito formula to:$\|Y_t-\overline{Y}_t\|^2$, we have:
\begin{eqnarray*}
d\|Y_t-\overline{Y}_t\|^2 &=&[\underbrace{2 \langle A_t-\overline{A}_t,Y_t-\overline{Y}_t \rangle }_1 \\ &+& \underbrace{2 \langle B_t Y_t-\overline{B}_t \overline{Y}_t,Y_t-\overline{Y}_t \rangle}_2\\ 
&+& \underbrace{2 \langle (C_t(Z_t)-\overline{C}_t(\overline{Z}_t)),Y_t-\overline{Y}_t \rangle }_3\\ &+&\|Z_t-\overline{Z}_t\|^2]dt+dM_t,   
\end{eqnarray*} 
where $dM_t:=2\langle Y_t-\overline{Y}_t,(Z_t-\overline{Z}_t)dW_t \rangle $.
\\By Cauchy-Schwarz inequality and Lipschitz continuity of $D_x\varphi$, \textbf{the first term 1} is greater than or equal to:
\begin{equation*}
-2 L^{x,\alpha}_{D_x\varphi} (\|X_t-\overline{X}_t\|+\| \alpha_t-\overline{\alpha}_t\|) \|Y_t-\overline{Y}_t\|.
\end{equation*}
\textbf{The second term 2} is equal to
\begin{equation*}
\begin{split}       
2\langle (B_t-\overline{B}_t)Y_t+\overline{B}_t(Y_t-\overline{Y}_t),Y_t-\overline{Y}_t \rangle.\end{split}
\end{equation*}
Since $Y$ is bounded, let $M_Y$ be its bound. By Cauchy-Schwarz inequality and the Lipschitz continuity condition of  $D_xb$, we have
\begin{equation*}
\begin{split}   
\langle (B_t-\overline{B}_t)Y_t,Y_t-\overline{Y}_t \rangle \geq -M_Y \normm{B_t-\overline{B}_t} \|Y_t-\overline{Y}_t\|,\\  \geq -M_Y L^{x,\alpha}_{D_xb} (\|X_t-\overline{X}_t\|+\| \alpha_t-\overline{\alpha}_t\|) \|Y_t-\overline{Y}_t\|,\end{split}
\end{equation*}
and by one-sided-Lipschitz condition of $b$, we have $$\langle \overline{B}_t(Y_t-\overline{Y}_t),Y_t-\overline{Y}_t \rangle \geq -c \|Y_t-\overline{Y}_t\|^2$$
So \textbf{the second term 2 } is greater than or equal to:
\begin{equation*}
\begin{split}
    -2M_Y L^{x,\alpha}_{D_xb}& (\|X_t-\overline{X}_t\|+\| \alpha_t-\overline{\alpha}_t\|) \|Y_t-\overline{Y}_t\|\\&-2c\|Y_t-\overline{Y}_t\|^2.
\end{split}
\end{equation*}

\item By the fact that $D_x\sigma$ does not depend on $x$, \textbf{the third term 3} is equal to  
\begin{equation*}    
2\langle C_t(Z_t)-C_t(\overline{Z}_t),Y_t-\overline{Y}_t \rangle.
\end{equation*}     
By linearity of the derivative and the trace, we have $$
C_t(Z_t)-C_t(\overline{Z}_t)=C_t(Z_t-\overline{Z}_t).
$$
And according to Lemma \RNum{4}.2, we have 
$$\|C_t(Z_t)-C_t(\overline{Z}_t)\| \leq L^{x}_{\sigma} \|Z_t-\overline{Z}_t\|.$$
By Cauchy-Schwarz, \textbf{the third term 3} is greater or equal to
\begin{equation*}
-(L^{x}_{\sigma})^2\|Y_t-\overline{Y}_t\|^2-\|Z_t-\overline{Z}_t\|^2.
\end{equation*}
\noindent We combine the previous results,  
we integrate between $s$ and $t$ ($s<t$), we take the expectation, we divide by $t-s$, and take the limit inferior as $s$ approaches to $t^-$, we obtain:
\begin{equation} \label{ineq}
D_-h^2(t)\geq 2\mu h^2(t)-2L_{Y}(h_x(t)+h_{\alpha}) h(t),
\end{equation}
where
\begin{equation*}
\begin{split}   
        h(t):=\|Y_t-\overline{Y}_t\|_{L^2},\text{ 
 } h_x(t):=\|X_t-\overline{X}_t\|_{L^2},\\h_{\alpha}:=\|\alpha-\overline{\alpha}\|_{\mathcal{A}},\text{    
 }L_Y:=L^{x,\alpha}_{D_x\varphi}+M_YL^{x,\alpha}_{D_xb}.
\end{split}
\end{equation*}
Now, if we put \begin{equation*}
f(t):=h(T-t),g(t):=L_Y(h_x(T-t)+h_{\alpha}),    
\end{equation*}
the inequality \eqref{ineq} becomes $$D^{+}f^{2}(t)\leq -2\mu f(t)^2+2f(t)g(t),$$
which gives after using \textbf{(Lemma III.4)}, and Gronwall's lemma (see the stability proof of the state process) $$f(t)\leq f(0) e^{-\mu t}+\int_0^t g(s) e^{-\mu(t-s)} ds.$$
After identifying the terms and using the fact that $D_x\psi$ is $L^{x}_{D_x\psi}$-Lipschitz continuous, we find the result. 
\end{dem}
\section{Contractivity of MSA}
\noindent Let's consider the function $\Tilde{H}:[0,T]\times\mathbb{R}^n\times\mathbb{R}^n\times A\rightarrow \mathbb{R}^{n}$ defined by
$$\Tilde{H}(t,x,y,a)=y^Tb(t,x,a)+\varphi(t,x,a).$$
since $\operatorname{Tr}(\sigma(t,x)^Tz)$ does not depend on $a$, so for all $t,x,y,z$ we have: $$\arg \min_{a\in A} H(t,x,y,z,a)=\arg \min_{a\in A} \Tilde{H}(t,x,y,a).$$ 
Suppose that for all $t, x, y$, the minimizer of the function \\$a \rightarrow \tilde{H}(t, x, y, a)$ (which exists since it's a continuous function on a compact set) is unique, and denotes $h(t, x, y)$.

\subsection{The continuity of h } 
\noindent

\begin{theorem}
\noindent \\The function $h:[0,T]\times \mathbb{R}^n\times \mathbb{R}^n\rightarrow A$ is continuous.
\end{theorem}
\begin{dem}\noindent
Let $\mathcal{X}:=[0,T]\times \mathbb{R}^n\times \mathbb{R}^n$, and let $(x_n)_{n\in \mathbb{N}}$ be a sequence of elements from $\mathcal{X}$ that converges to $x \in \mathcal{X}$. We have that $$h(x_n)\xrightarrow[n\to \infty]{} h(x).$$
We only need to prove that every subsequence of $(h(x_n))_{n\in \mathbb{N}}$ has a convergent subsequence that converges to $h(x)$.
\noindent \\ Let $(h(x_{\phi_1(n)}))_{n\in \mathbb{N}}$ be a subsequence of $(h(x_n))_{n\in \mathbb{N}}$. It is a sequence of elements of $A$ which is a compact, so it has a convergent subsequence $(h(x_{\phi_1(\phi_2(n)}))_{n\in \mathbb{N}}$ let $l\in A$ be its limit. We have, by definition of $h$: $$\forall a\in A, \Tilde{H}(x_{\phi_1(\phi_2(n))},h(x_{\phi_1(\phi_2(n))})\leq \Tilde{H}(x_{\phi_1(\phi_2(n))},a).$$
By taking the limit, we obtain $$\forall a\in A, \Tilde{H}(x,l)\leq \Tilde{H}(x,a).$$
Due to the uniqueness of the minimizer of $\Tilde{H}(x,.)$, we have $l=h(x)$. Hence, the continuity of $h$.
\end{dem}
\subsection{The Lipschitz continuity of MSA }
\noindent
Note that $\operatorname{MSA}:\mathcal{A}\rightarrow\mathcal{A}$ is well defined by continuity of $h$.
\begin{assumption}
Now assume that $h$ is Lipschitz continuous with respect to $(x,y)$ uniformly in $t$, and let $L_h$ be its constant.
\end{assumption}
\begin{theorem}
MSA is Lipschitz continuous. $$\|\operatorname{MSA}(\alpha^2)-\operatorname{MSA}(\alpha^1)\|_{\mathcal{A}}\leq L_{\mu,T} \|\alpha^2-\alpha^1\|_{\mathcal{A}},$$
where 
\begin{equation*}
L_{\mu,T}:=L_h\left[(L^{\alpha}_{b}(1+L^{x}_{D_x\psi})+L_Y) M_{\mu,T}+L_YL^{\alpha}_{b}M_{\mu,T}^2\right]
\end{equation*}
where $M_{\mu,T}:=\displaystyle\frac{1-e^{-\mu T}}{\mu}$
\end{theorem}
\begin{dem}
\noindent
\noindent 
Let $\alpha^1,\alpha^2\in \mathcal{A}$, $X^1,X^2$ the two associated state processes and $Y^1,Y^2$ the two associated adjoint processes
\begin{equation*}
\begin{split}        
\|\operatorname{MSA}&(\alpha^2)-\operatorname{MSA}(\alpha^1)\|_{\mathcal{A}}=\\&\sup_t\sqrt{\mathbb{E}[\|h(t,X^2_t,Y^2_t)-h(t,X^1_t,Y^1_t)\|^2]} ,  
\end{split}
\end{equation*}
since $h$ is Lipschitz continuous, we have $$\|h(t,X^2_t,Y^2_t)-h(t,X^1_t,Y^1_t)\|\leq L_h(\|X^2_t-X^1_t\|+\|Y^2_t-Y^1_t\|).$$
Then $$\|\operatorname{MSA}(\alpha^2)-\operatorname{MSA}(\alpha^1)\|_{\mathcal{A}}\leq L_h(\|X^2-X^1\|_{\mathcal{A}}+\|Y^2-Y^1\|_{\mathcal{A}}).$$
According to \eqref{state stability}, we have 
\begin{equation*} 
\begin{split}\|X^2_t-X^1_t\|_{L^2}\leq & \|X^2_0-X^1_0\|_{L^2}e^{-\mu t}\\&+L^{\alpha}_{b}\frac{1-e^{-\mu t}}{\mu} \|\alpha^2-\alpha^1\|_{\mathcal{A}}.
\end{split}
\end{equation*}
The initial condition of the state process is not changed in the MSA algorithm. Therefore $X^2_0=X^1_0$, and 
$$\|X^2-X^1\|_{\mathcal{A}}\leq L^{\alpha}_{b}M_{\mu,T} \|\alpha^2-\alpha^1\|_{\mathcal{A}}.$$
According to \eqref{adjoint stability}, we have 
\begin{equation*}
\begin{split}    
&\|Y^2_t-Y^1_t\|_{L^{2}}\leq L^{x}_{D_x\psi} \|X^2_T-X^1_T\|_{L^{2}} e^{-\mu (T-t)}+\\&L_Y \displaystyle\int_t^T (\|X^2_{T+t-s}-X^1_{T+t-s}\|_{L^{2}}+\|\alpha^2-\alpha^1\|_{\mathcal{A}}) e^{-\mu (T-s)}ds.
\end{split}
\end{equation*}
Therefore,
\begin{equation*}
\begin{split}
\|Y^2-Y^1\|_{\mathcal{A}}\leq&\left[(L^{\alpha}_{b}L^{x}_{D_x\psi}+L_Y)M_{\mu,T}+L_YL^{\alpha}_{b}M_{\mu,T}^2\right]\\&\|\alpha^2-\alpha^1\|_{\mathcal{A}} .   
\end{split}\end{equation*}
and $$\|\operatorname{MSA}(\alpha^2)-\operatorname{MSA}(\alpha^1)\|_{\mathcal{A}}\leq L_{\mu,T}\|\alpha^2-\alpha^1\|_{\mathcal{A}}.$$
\end{dem}
\noindent
\begin{theorem}
MSA is a contraction:
\begin{enumerate}
    \item For large values of $\mu$ no matter $T$, or
    \item For small values of $T$ no matter $\mu$.
    \end{enumerate}
\end{theorem}
\begin{dem}\noindent
\begin{enumerate}
    \item Let $\overline{T}\in \mathbb{R}_{\geq 0}\cup \{+\infty\}$. We have $0 \leq  M_{\mu,T}\leq \frac{1}{\mu}$, so $M_{\mu,T}\to 0$ as $( \mu,T)\to (+\infty,\overline{T} )$ so $$L_{\mu,T} \xrightarrow[( \mu,T )\to (+\infty,\overline{T})]{} 0$$ 
    \item Let $ \overline{ \mu}\in \mathbb{R}_{\geq 0}\cup \{ +\infty\}$. The case $\overline { \mu} = + \infty$ has already been dealt with in 1. Therefore, assume that $\overline{\mu}\in \mathbb{R}_{\geq 0}$ we have $M_{\mu,T}\underset{(\overline{\mu},0^+)}{\sim} T$ so $M_{\mu,T}\to 0$ when $(\mu,T)\to (\overline{\mu}^{+},0)$ 
    so $$L_{\mu,T}\xrightarrow[(\mu,T)\to (\overline{\mu}^{+},0)]{} 0$$
\end{enumerate}
\end{dem}
\subsection{The convergence of MSA }
\noindent 
\noindent
\begin{theorem}
Let $\mu$ and $T$ such that $L_{\mu,T}<1$.\\
The MSA algorithm converges, and we have:
$$\|\alpha^n-\alpha\|_{\mathcal{A}}\leq \frac{(L_{\mu,T})^n}{1-L_{\mu,T}} \|\alpha^1-\alpha^0\|_{\mathcal{A}},$$
where $\alpha$ is the optimal control.
\end{theorem}
\begin{dem}
MSA is contracting in $\mathcal{A}$ which is a closed set of a Banach space (the set of adapted processes $\alpha$ such that $\|\alpha\|_{\mathcal{A}}$ is finite). Therefore, according to fixed point theorem of Banach-Picard, it has a unique fixed point which is also the limit of the sequence $(\alpha^{n})_{n \in \mathbb{ N}}$.

\end{dem}
\bibliographystyle{IEEEtran}
\bibliography{Ref}
\end{document}